\documentclass[12pt, reqno]{amsart}
\usepackage{amsmath,amsthm,amsfonts,amssymb,enumerate,amscd}
\newtheorem{Defin}{Definition}[section]
 \newtheorem{Theo}[Defin]{Theorem}
 \newtheorem{Lem}[Defin]{Lemma}
 \newtheorem{Cor}[Defin]{Corollary}
 
 \newtheorem{Rem}[Defin]{Remark}
 \newtheorem{Pro}[Defin]{Proposition}
 \newtheorem{Def}[Defin]{Definition}
 \newtheorem{Conj}[Defin]{Conjecture}
 \newenvironment{demo}{\textsl{Proof. }\rm}{$\hfill\square$}
 
 \renewcommand{\thefootnote}{}

\DeclareMathOperator{\Div}{\rm div}

\def\bbr{\mathbb R}

\def\calc{\mathcal C}

\def\rme{\mathrm e}

\def\Saut{\vskip 0.5 true cm \noindent}

\usepackage[dvips]{graphicx}
 \numberwithin{equation}{section} \voffset -1in \hoffset -1in
 \catcode`\@=11
\catcode`\@=\active
\def \logo@{\relax}
\usepackage{color}

\newcommand{\rd}{\textcolor{red}}

\def\*{{\phantom *}}
\begin{document}
\phantom{.} {\qquad \hfill \rd{\textit{\textbf{\textit{Lyapunov
Memorial Conference}, June 24-30, 2007}}}}

\vskip 1.0cm

\setcounter{footnote}{0}
\renewcommand{\thefootnote}{\arabic{footnote}}

\title[Dirichlet-to-Neumann (Gibbs) Semigroups ]
{From Laplacian Transport to Dirichlet-to-Neumann (Gibbs) Semigroups}
\maketitle

\Saut

\begin{center}

\vspace{0.4cm}

\textbf{Valentin A.Zagrebnov} \footnote{E-mail: Valentin.Zagrebnov$@$cpt.univ-mrs.fr}

\vspace{0.2cm}

\textbf{Universit\'e de la M\'editerran\'ee \\ Centre de Physique
Th\'eorique - UMR 6207 \footnote{Universit\'{e} de Provence - Aix-Marseille I,
Universit\'{e} de la M\'{e}diterran\'{e}e - Aix-Marseille II,
Universit\'{e} du Sud - Toulon - Var, FRUMAM (FR 2291)} \\ Luminy - Case 907,
13288 Marseille, Cedex 09, France}

\vspace{0.5cm}

\textbf{Abstract}\footnote{ An extended version of the author's talk presented on the \textit{Lyapunov
Memorial Conference}, June 24-30, 2007 (Kharkov University, Ukraine), which is based on the common project
with Professor Hassan Emamirad (Laboratoire de Math\'{e}matiques, Universit\'{e} de Poitiers).}

\end{center}

\noindent The paper gives a short account of some basic properties of \textit{Dirichlet-to-Neumann}
operators $\Lambda_{\gamma,\partial\Omega}$ including the corresponding semigroups motivated by the Laplacian
transport in anisotropic media ($\gamma \neq I$) and by elliptic systems with dynamical boundary conditions.
For illustration of these notions and the properties we use the explicitly constructed \textit{Lax semigroups}.
We demonstrate that for a general smooth bounded convex domain $\Omega \subset \mathbb{R}^d$ the corresponding
{Dirichlet-to-Neumann} semigroup $\left\{U(t):= e^{-t \Lambda_{\gamma,\partial\Omega}}\right\}_{t\geq0}$ in the
Hilbert space $L^2(\partial \Omega)$  belongs to the \textit{trace-norm} von Neumann-Schatten ideal for any $t>0$.
This means that it is in fact an \textit{immediate Gibbs} semigroup. Recently Emamirad and Laadnani have
constructed a \textit{Trotter-Kato-Chernoff} product-type approximating family
$\left\{(V_{\gamma, \partial\Omega}(t/n))^n \right\}_{n \geq 1}$ \textit{strongly} converging to the semigroup
$U(t)$ for $n\rightarrow\infty$. We conclude the paper by discussion of a conjecture about convergence of the
\textit{Emamirad-Laadnani approximantes} in the the {\textit{trace-norm}} topology.

\vspace{0.2cm}


\newpage
\tableofcontents

\section{Laplacian transport and Dirichlet-to-Neumann operators}\label{Intro}

\smallskip

\noindent \textbf{Example 1.1.} Is is well-known (see e.g. \cite{LeUh}) that the problem of determining
a \textit{conductivity matrix} field $\gamma(x) = [\gamma_{i,j}(x)]_{i,j=1}^{d}$, for $x$ in a bounded
open domain $\Omega \subset \mathbb{R}^d$, is related to "measuring" the elliptic \textit{Dirichlet-to-Neumann} map
for associated conductivity equation. Notice that solution of this problem has a lot of practical applications in
various domains: geophysics, electrochemistry etc. It is also an important diagnostic tool in medicine, e.g. in the
\textit{electrical impedance tomography}; the tissue in the human body is an example of highly anisotropic
conductor \cite{BaBr}.

Under assumption that there is no sources or sinks of current the potential $v(x), \ x\in \Omega,$ for a
given voltage $f(\omega), \ \omega\in \partial\Omega,$ on the (smooth) boundary $\partial\Omega$ of $\Omega$ is
a solution of the Dirichlet problem:
$$\begin{cases}
 \Div(\gamma\nabla v) = 0 \quad \text{in}\;\;\Omega,\\
 v|_{\partial\Omega} = f\quad\text{on}\;\; \partial \Omega.\end{cases}
 \leqno{\mathrm{\bf (P1)}}$$
Then the corresponding to (\textbf{P1}) Dirichlet-to-Neumann map (operator) $\Lambda_{\gamma,\partial\Omega}$
is defined by
\begin{equation}\label{DN-oper-Lambda}
\Lambda_{\gamma,\partial\Omega}: f \mapsto  \partial v_f/\partial
\nu_\gamma := \nu\cdot\gamma \ \nabla v_f\mid_{\partial\Omega} \ .
\end{equation}
Here $\nu$ is the unit outer-normal vector to the boundary at $\omega \in\partial\Omega$
and the function $u :=u_f$ is solution of the Dirichlet problem (\textbf{P1}).

The Dirichlet-to-Neumann operator (\ref{DN-oper-Lambda}) is also called the \textit{voltage-to-current}
map, since the function $\Lambda_{\gamma,\partial\Omega} f$ gives the induced current flux trough the boundary
$\partial\Omega$. The key (\textit{inverse}) problem is whether on can determine the conductivity matrix $\gamma$
by knowing electrical boundary measurements, i.e. the corresponding Dirichlet-to-Neumann operator? Unfortunately,
this operator does not determine the matrix $\gamma$ uniquely, see e.g. \cite{GrUh} and references there.

\bigskip

\noindent \textbf{Example 1.2.} The problem of electrical current flux in the form (\textbf{P1}) is an
example of so-called \textit{Laplacian transport}. Besides the voltage-to-current problem the motivation
to study this kind of transport comes for instance from the \textit{transfer} across biological membranes,
see e.g. \cite{Sap}, \cite{GrFiSap}.

Let some "species" of concentration $C(x)$, $x\in\mathbb{R}^d$, diffuse in the
\textit{isotropic} bulk ($\gamma = I$) from a (distant) source localised on the closed boundary
$\partial\Omega_0$ towards a \textit{semipermeable} compact interface $\partial\Omega$ on
which they disappear at a given rate $W$. Then the \textit{steady} concentration field (Laplacian
transport with a diffusion coefficient $D$) obeys the set of equations:
$$\begin{cases}
\Delta C = 0 , \  x \in\Omega_0 \setminus \Omega \ , \\
C(\omega_0\in \partial\Omega_0) = C_0 , \ {\rm{at \ the \ source}} \ , \\
D \ \partial_\nu C(\omega) = W \ (C(\omega) - 0), \ {\rm{on \ the \ interface}} \ \omega\in \partial\Omega \
\end{cases}
\leqno{\mathrm{\bf (P2)}}$$
Let $C=C_0 (1-u)$. Then $\Delta u$ = 0, $x\in \Omega$. If we put $\mu:= D/W$, then the boundary conditions
on $\partial\Omega$ take the form:
$(I - \mu \partial_\nu)u\mid_{\partial\Omega}(\omega)=1\mid_{\partial\Omega}(\omega)$, where
$(1 \mid_{\partial\Omega})(\omega) = \chi_{\partial\Omega} (\omega)$ is characteristic function of
the set $\partial\Omega$, and  $u(\omega_0) = 0, \omega_0\in \partial\Omega_0$ on the source boundary.

Consider now the following auxiliary Laplace-Dirichlet problem:
\begin{equation}\label{L-D-Probl}
\Delta u = 0 , \  x \in\Omega_0 \setminus \Omega \ ,  \  \
u \mid_{\partial\Omega} (\omega)= f (\omega \in\partial\Omega) \ \
{\rm{and}} \ \ u \mid_{\partial\Omega_0} (\omega)= 0 \ ,
\end{equation}
with solution $u_f$. Then similar to (\ref{DN-oper-Lambda}) we can associate with the problem
(\ref{L-D-Probl}) a Dirihlet-to-Neumann operator
\begin{equation}\label{D-to-N-Lapl-transp}
\Lambda_{\gamma=I,\partial\Omega}: f \mapsto \partial_\nu u_f \mid_{\partial\Omega}
\end{equation}
with domain $\mathrm{dom}(\Lambda_{I,\partial\Omega})$, which belongs to a certain \textit{Sobolev} space, 
Section \ref{D-to-N}.

The advantage of this approach is that as soon as the operator (\ref{D-to-N-Lapl-transp}) is defined one
can apply it to study the \textit{mixed} boundary value problem (\textbf{P2}). This gives in particular the
value of the particle flux due to Laplacian transport across the membrane $\partial\Omega$.
Indeed, one obtain that  $(I + \mu \Lambda_{I,\partial\Omega})u \mid_{\partial\Omega} =
1 \mid_{\partial\Omega}$,
and that the local (diffusive) particle {\textit{flux}} is defined as:
\begin{equation}\label{flux}
\phi \mid_{\partial\Omega}:=
D \ C_0(-\partial_n u)\mid_{\partial\Omega}=
D \ C_0(\Lambda_{I,\partial\Omega}(I+\mu\Lambda_{I,\partial\Omega})^{-1}1) \mid_{\partial\Omega} \ .
\end{equation}
Then the corresponding total flux across the membrane $\partial\Omega$ :
\begin{equation}\label{total-flux}
\Phi:= (\phi, 1)_{L^2 (\partial\Omega)}=
D \ C_0(\Lambda(I+\mu\Lambda_{I,\partial\Omega})^{-1}1 , 1)_{L^2 (\partial\Omega)}
\end{equation}
is experimentally measurable macroscopic response of the system, expressed via
transport parameters $D, C_0, \mu$ and geometry of $\partial\Omega$. Here
$(\cdot, \cdot)_{L^2 (\partial\Omega)}$
is scalar product in the Hilbert space $\partial \mathcal{H}:= L^2 (\partial\Omega)$.

The aim of the present paper is twofold: \\
(i) to give a short account of some standard results
about  Dirichlet-to-Neumann operators and related \textit{Dirichlet-to-Neumann semigroups} that solve
a certain class of elliptic systems with dynamical boundary conditions;\\
(ii) to present some recent results concerning the \textit{approximation} theory and the \textit{Gibbs}
character of the Dirichlet-to-Neumann semigroups for compact sets $\Omega$ with smooth
boundaries $\partial\Omega$.

To this end in the next Section \ref{D-to-N} we recall some fundamental properties of the
Dirichlet-to-Neumann operators and semigroups, we illustrate them by few elementary examples,
including the \textit{Lax semigroups} \cite{Lax}.

In Section \ref{Appr} we present the strong \textit{Emamirad-Laadnani approximations} of the Dirichlet-to-Neumann
semigroups inspired by the \textit{Chernoff} theory and by its generalizations in \cite{NeZag}, \cite{CaZag2}.

We show in Section \ref{Gibbs} that for compact sets $\Omega$ with smooth boundaries $\partial\Omega$
the Dirichlet-to-Neumann semigroups are in fact (immediate) \textit{Gibbs} semigroups \cite{Zag2}.

Some recent results and conjectures about approximations of the Dirichlet-to-Neumann (Gibbs) semigroups
in operator and \textit{trace-norm} topologies are collected in the last Section \ref{Concl}.


\section{Dirichlet-to-Neumann operators and semigroups}\label{D-to-N}


\noindent \textbf{2.1 Dirichlet-to-Neumann operators}

\smallskip

Let  $\Omega$ be an open bounded domain in $ \mathbb R^d$ with a smooth boundary $\partial \Omega$.
Let $\gamma$ be a $C^\infty(\overline \Omega)$ matrix-valued function on $\overline \Omega$, which
we call the \textit{Laplacian transport matrix} in domain $\Omega$.

We suppose that the matrix-valued function $\gamma(x):= [\gamma_{i,j}(x)]_{i,j=1}^{d}$ satisfies the
following hypotheses:\\
\textbf{(H1)} The real coefficients are symmetric and $\gamma_{i,j}
(x)=\gamma_{j,i}(x)\in \calc^\infty(\overline \Omega)$.\\
\textbf{(H2)} There exist two constants $0 < c_1 \le c_2 < \infty$ such that for all
$\xi \in \mathbb R^d$, we have
\begin{equation}\label{minmax}
c_1\|\xi\|^2 \le \sum_{i,j=1}^n \xi_i\xi_j\gamma_{i,j}(x)\le c_2\|\xi\|^2.
\end{equation}

Then the \textit{Dirichlet-to-Neumann} operator $\Lambda_{\gamma,\partial\Omega} $ associated with
the Laplacian transport in $\Omega$ is defined as follows.

Let $f\in C(\partial \Omega)$ and denote by $v_f$
the \textit{unique} solution (see e.g. \cite{GiTr}, Theorem 6.25) of the Dirichlet problem
$$\begin{cases}
 A_{\gamma,\partial\Omega}\ v:=\Div(\gamma \ \nabla v) = 0 \quad \text{in}\;\;\Omega \ ,\\
 v\mid_{\partial\Omega} = f\quad\text{on}\;\; \partial \Omega \ ,\end{cases}
 \leqno{\mathrm{\bf (P1)}}\label{P}$$
in the Banach space $X:= C(\overline \Omega)$. Here operator $A_{\gamma,\partial\Omega} $ is defined on its
maximal domain
\begin{equation}\label{dom-A-X}
{\rm{dom}}(A_{\gamma,\partial\Omega}): = \{u\in X : A_{\gamma,\partial\Omega} \ u \in X\} \ .
\end{equation}
\begin{Def}\label{DtoN-X}
The Dirichlet-to-Neumann operator is the map{\rm{:}}
\begin{equation}\label{Lambda}
\Lambda_{\gamma,\partial\Omega}: f \mapsto  \partial v_f/\partial
\nu_\gamma= \nu\cdot\gamma\nabla v_f\mid_{\partial\Omega} \ ,
\end{equation}
with domain {\rm{:}}
\begin{eqnarray}\label{dom-Lambda}
&& {\rm{dom}}(\Lambda_{\gamma,\partial\Omega})= \\
&& \{f\in \partial C(\Omega_R): v_f\in{\rm{Ker}}(A_{\gamma,\partial\Omega}) \ \ {\rm{and}} \ \
|(\nu\cdot\gamma\nabla v_f\mid_{\partial\Omega})| < \infty \} . \nonumber
\end{eqnarray}
Here $\nu$ denotes the unit outer-normal vector at $\omega\in\partial\Omega$
and  $v_f$ is the solution of Dirichlet problem {\rm{\textbf{(P1)}}}.
\end{Def}
The solution $v_f: = L_{\partial\Omega} f$ of the problem
\textbf{(P1)} is called the $\gamma$-\textit{harmonic lifting} of $f$, where
$L_{\partial\Omega}: C(\partial \Omega)\mapsto C^2 (\Omega)\cap C(\overline \Omega)$ is called the
\textit{lifting} operator with domain ${\rm{dom}}(L_{\partial\Omega})= C(\partial \Omega)$.
If $T_{\partial\Omega} : C(\overline \Omega)\mapsto C(\partial \Omega)$ denotes the \textit{trace}
operator on
the smooth boundary $\partial \Omega$, i.e. $v\mid_{\partial\Omega} = T_{\partial\Omega} \ v$, then
\cite{Eng}:
\begin{equation}\label{L-T-oper}
L_{\partial\Omega}= (T_{\partial\Omega}\mid_{{\rm{Ker}}(A_{\gamma,\partial\Omega})} \, )^{-1} \
\ {\rm{and}} \ \
{\rm{dom}}(\Lambda_{\gamma,\partial\Omega})= T_{\partial\Omega}\{{\rm{Ker}}(A_{\gamma,\partial\Omega})\} \ .
\end{equation}
\begin{Rem}\label{L-T-oper-Rem}
Let $\partial X := C(\partial \Omega)$. Then (\ref{L-T-oper}) implies :
\begin{equation}\label{Ident}
T_{\partial\Omega} L_{\partial\Omega} \ u = u \ , \ u \in \partial X  \ \
{\rm{and}} \ \ L_{\partial\Omega} T_{\partial\Omega}\ w = w \ , \ w\in {\rm{Ker}}
(A_{\gamma,\partial\Omega}) \ .
\end{equation}
One also gets that the lifting operator is bounded:  $L_{\partial\Omega}\in \mathcal{L}(\partial X, X)$,
whereas
the Dirichlet-to-Neumann operator (\ref{Lambda}) is obviously not.
\end{Rem}

Now let $\mathcal{H}$ be Hilbert space $L^2(\Omega)$ and $\partial
\mathcal{H}:=L^2(\partial \Omega)$ denote the {\it boundary space}.
In order that the problem {\bf (P1)} admits a \textit{unique} solution $v_f$,
one has to assume that $f\in W^{1/2}_{2}(\partial \Omega)$, and then $v_f$ belongs the \textit{Sobolev} space
$W^{1}_{2}(\Omega)$, see e.g. \cite[Ch.7]{Tay}. So, we can define Dirichlet-to-Neumann operator in
the Hilbert space $\partial \mathcal{H}$ by (\ref{Lambda}) with domain:
\begin{equation}\label{DtoN-H-dom}
{\rm{dom}}(\Lambda_{\gamma,\partial\Omega}):= \{ f\in W^{1/2}_{2}(\partial \Omega): \
\Lambda_{\gamma,\partial\Omega} f \in \partial \mathcal{H}= L^2(\partial \Omega)\}.
\end{equation}
\begin{Pro}\label{DtoN-H-Prop}
The Dirichlet-to-Neumann operator (\ref{Lambda}) with domain (\ref{DtoN-H-dom}) in the Hilbert space
$\partial \mathcal{H}$ is unbounded, non-negative, self-adjoint, first-order elliptic pseudo-differential
operator with compact resolvent.
\end{Pro}
\noindent The complete proof can be found e.g. in \cite[Ch.7]{Tay}, \cite{Tay1}. Therefore, we give
here only some comments on these properties of the Dirichlet-to-Neumann operator (\ref{Lambda})
in $\partial \mathcal{H}= L^2(\partial \Omega)$.
\begin{Rem}\label{DtoN-H-Properties}
{\rm{(a)}} By virtue of definition (\ref{Lambda}) for any $f\in W^{1/2}_{2}(\partial \Omega)$ one gets:
\begin{eqnarray}\label{Posit}
&& (f , \Lambda_{\gamma,\partial\Omega} f)_{\partial \mathcal{H}} =
\int_{\partial \Omega}\ d\sigma(\omega)\ \overline{v_f(\omega)} \ \nu\cdot\gamma(\omega)(\nabla v_f)(\omega)
= \\ \nonumber
&& \int_{\Omega}\ dx \ \Div(\overline{v_f(x)} \ (\gamma \nabla v_f)(x)) =
\int_{\Omega}\ dx \ (\nabla \overline{v_f(x)} \cdot \gamma \ \nabla v_f)(x)) \geq 0 \ ,
\end{eqnarray}
since the matrix $\gamma$ verifies {\rm{\textbf{(H2)}}}. Thus, operator $\Lambda_{\gamma,\partial\Omega}$ is
non-negative.\\
{\rm{(b)}} In fact to ensure the existence of the trace
$T_{\partial\Omega} (\nu\cdot\gamma \nabla (L_{\partial\Omega} f))$
one has initially to define operator $\Lambda_{\gamma,\partial\Omega}$ for
$f\in W^{3/2}_{2}(\partial \Omega)$.
Then Dirichlet-to-Neumann operator is a self-adjoint extension with domain (\ref{DtoN-H-dom})
and moreover it is a bounded map
$\Lambda_{\gamma,\partial\Omega}: W^{1/2}_{2}(\partial \Omega)\mapsto W^{-1/2}_{2}(\partial \Omega)$. \\
{\rm{(c)}} By (\ref{Posit}) and since derivatives of the first-order are involved in (\ref{Lambda}) one
can conclude that this operator should be elliptic and pseudo-differential. If $\gamma(x) = I$, then
$\Lambda_{I,\partial\Omega}$ is, roughly, the operator $(-\Delta_{\partial\Omega})^{1/2}$,
where $\Delta_{\partial\Omega}$ is the Laplace-Beltrami operator on $\partial\Omega$, with corresponding
induced metric \cite[Ch.7]{Tay}, \cite{Tay1}.\\
{\rm{(d)}} Compactness of the imbedding $W^{1/2}_{2}(\partial \Omega) \hookrightarrow
L^2(\partial \Omega)$ implies the compactness of the resolvent of  $\Lambda_{\gamma,\partial\Omega}$.
\end{Rem}
By (a) and (d) the spectrum $\sigma(\Lambda_{\gamma,\partial\Omega})$ of the Dirichlet-to-Neumann operator
is a set of non-negative increasing eigenvalues $\{\lambda_k\}_{k=1}^{\infty}$. The rate of increasing is
given by the Weyl asymptotic formula, see e.g. \cite{Hor}, \cite{Tay}:
\begin{Pro}\label{Weyl-asympt}
Let $\Lambda_{\gamma,\partial\Omega}(x,\xi)$, for $(x,\xi)\in T^*\partial\Omega$, be the
symbol of the first-order, elliptic pseudo-differential Dirichlet-to-Neumann operator
$\Lambda_{\gamma,\partial\Omega}$. Then the asymptotic behaviour of the
corresponding eigenvalues as $k\to \infty$ has the form:
\begin{equation*}
\lambda_k \sim \left\{\frac k{C(\partial \Omega,\Lambda_\gamma)}\right\}^{1/(d-1)} \ ,
\end{equation*}
where
\begin{equation*}
C(\partial \Omega,\Lambda_\gamma) : = \frac
1{(2\pi)^{d-1}}\int_{\Lambda_{\gamma,\partial\Omega}(x,\xi)\le 1} dx\;d\xi \ .
\end{equation*}
\end{Pro}
Another important result is due to Hislop and Lutzer \cite{HiLu}. It concerns a localization (\textit{rapid decay})
of the $\gamma$-harmonic lifting of the corresponding eigenfunctions.
\begin{Pro}\label{Localization}
Let $\{\phi_k\}_{k=1}^{\infty}$ be eigenfunctions of the Dirichlet-to-Neumann operator:
$\Lambda_{\gamma,\partial\Omega}= \lambda_k\phi_k$ with $\|\phi_k\|_{L^2(\partial \Omega)} = 1$.
Let $v_{\phi_k}: = L_{\partial\Omega} \phi_k$ be $\gamma$-harmonic lifting of $\phi_k$ to $\Omega$
corresponding to the problem {\rm{\textbf{(P1)}}}.
Then for any compact $\mathcal{C} \subset \phi_k$ and $x\in \mathcal{C}$ one gets the representation:
\begin{equation}\label{estim}
|v_{\phi_k}(x)| =
\psi(x,p,\mathcal{C}) / \lambda_k \, ^{p}
\end{equation}
with arbitrary large $p >0$. Here $\psi(x,p,\mathcal{C})$ is a decreasing function of the distance
$\mathrm{dist}(x ,\, \partial\Omega)$.
\end{Pro}
\noindent Since by the Weyl asymptotic formula we have $\lambda_k = O(k^{1/(d-1)})$, the decay implied by
the estimate (\ref{estim}) is algebraic.
\begin{Conj}\label{Exp-loc} \cite{HiLu} In fact the order of decay instead of
$\psi(x,p,\mathcal{C}) / \lambda_k \, ^{p}$ is exponential:
$O(\exp [- \, k \, \mathrm{dist}(\mathcal{C} ,\, \partial\Omega)])$ .
\end{Conj}

\bigskip
\noindent \textbf{2.2 Example of a Dirichlet-to-Neumann operator}
\smallskip

To illustrate the results mentioned above we consider a simple example which will be useful below for
contraction of the \textit{Lax semigroups}.

Consider a homogeneous isotropic case: $\gamma(x)= I$, and let $\Omega = \Omega_R : =
\{x \in \mathbb{R}^{d=3}: \ \|x\| < R\}$.
Then $A_{\gamma,\partial\Omega_R} = \Delta_{\partial\Omega_R}$ and for the \textit{harmonic lifting} of
\begin{equation*}
f(\omega)= \sum_{l,m} f_{l,m}^{(R)}  \ Y_{l,m} (\theta, \varphi) \in W^{1/2}_{2}(\partial \Omega_R) \ ,
\end{equation*}
we obtain:
\begin{equation}\label{v_f=R}
v_f(r, \theta, \varphi)= \sum_{l,m} \left(\frac{r}{R}\right)^l \ f_{l,m}^{(R)}  \ Y_{l,m} (\theta, \varphi),
\end{equation}
since the \textit{spherical} functions $ \ \{Y_{l,m}\}_{l=0, |m|\leq l}^{\infty} \ $ form a
complete orthonormal
basis in the Hilbert space $\partial \mathcal{H}= L^2(\partial \Omega_R,
\ d\theta \, \sin\theta \, d\varphi)$.

Definition (\ref{Lambda}) and (\ref{v_f=R}) imply that non-negative, self-adjoint, first-order elliptic
pseudo-differential Dirichlet-to-Neumann operator
\begin{equation}\label{DtoN=R}
(\Lambda_{I,\partial\Omega_R} f) (\omega= (R, \theta, \varphi)) =
\sum_{l=0}^{\infty}\sum_{m=-l}^{m=l} \left(\frac{l}{R}\right) \ f_{l,m}^{(R)}  \ Y_{l,m} (\theta, \varphi) ,
\end{equation}
has discrete spectrum $\sigma(\Lambda_{I,\partial\Omega_R}) := \{\lambda_{l,m}= \ l/R\}_{l=0,
|m|\leq l}^{\infty}$
with spherical eigenfunctions:
\begin{equation}\label{DtoN=R-spect}
(\Lambda_{I,\partial\Omega_R} Y_{l,m}) (R, \theta, \varphi) = \left(\frac{l}{R}\right)\
Y_{l,m}(\theta, \varphi) ,
\end{equation}
and multiplicity $m$. The operator (\ref{DtoN=R}) is obviously unbounded and it has a compact resolvent.
\begin{Rem}\label{DtoN-H-localization}
Since by virtue of (\ref{v_f=R}) the $\gamma$-harmonic lifting of the eigenfunction $Y_{l,m}$ to the ball
$\Omega_R$ is
\begin{equation*}\label{R-lifting}
v_{\, Y_{l,m}}(r, \theta, \varphi)= \left(\frac{r}{R}\right)^l \ Y_{l,m} (\theta, \varphi) \ ,
\end{equation*}
one can check the localization (Proposition \ref{Localization}) and Conjecture about the exponential
decay explicitly. For distances: $0 < \mathrm{dist}(x, \partial\Omega_R)= R - r \ll R $, one obtains
$|v_{\, Y_{l,m}}(r, \theta, \varphi)| = O(e^{- l (R - r)/R})$.
\end{Rem}
\smallskip
\noindent \textbf{2.3 Dirichlet-to-Neumann semigroups on $\partial X$}
\smallskip

To define the Dirichlet-to-Neumann semigroups on the \textit{boundary} Banach space
$\partial X = C(\partial\Omega)$ we can follow  the line of reasoning of \cite{Esc} or \cite{Eng}.
To this end consider in $X = C(\Omega)$ the following elliptic system with \textit{dynamical boundary conditions}
$$\begin{cases}
\Div(\gamma\nabla u(t,\cdot))= 0
\quad &\text{in}\;\;(0,\infty)\times\Omega,\\
\partial u(t,\cdot)/{\partial t} + \partial u(t,\cdot)/{\partial
\nu_\gamma}=0 &\text{on}\;\;(0,\infty)\times\partial\Omega,\\
u(0,\cdot) = f\quad&\text{on}\;\; \partial \Omega.\end{cases}
\leqno{\mathrm{\bf (P2)}}\label{P2}$$
\begin{Pro}\label{DtoN-X-semigr-def}
The problem {\rm{\textbf{(P2)}}} has a unique solution $u_f(t,x)$ for any $f\in C(\partial\Omega)$. Its
trace on the boundary $\partial\Omega$ has the form:
\begin{equation}\label{solution}
u_f(t,\omega):= (T_{\partial\Omega}u_f(t,\cdot))(\omega) = (U(t) f)(\omega) \ ,
\end{equation}
where the family of operators $\{U(t) = e^{- t \Lambda_{\gamma,\partial\Omega}}\}_{t\geq0}$ is a
$C_0$-semigroup
generated by the Dirichlet-to-Neumann operator of the problem {\rm{\textbf{(P1)}}}.
\end{Pro}
The following key result about the properties of the Dirichlet-to-Neumann semigroups on the boundary
Banach space $\partial X = C(\partial\Omega)$ is due to Escher-Engel \cite{Esc},\cite{Eng}
and Emamirad-Laadnani \cite{EmLa}:
\begin{Pro}\label{DtoN-X-semigr-prop-PROP}
The semigroup $\{U(t) = e^{- t \Lambda_{\gamma,\partial\Omega}}\}_{t\geq0}$ is analytic, compact, positive,
irreducible and Markov $C_0$-semigroup of contractions on $C(\partial\Omega)$.
\end{Pro}
\begin{Rem}\label{DtoN-X-semigr-prop-REMARK}
The complete proof can be found in the papers quoted above. So, here we make only some comments and hints
concerning the Proposition \ref{DtoN-X-semigr-prop-PROP}.
\end{Rem}

\smallskip
\noindent \textbf{2.4 Dirichlet-to-Neumann semigroups on $\partial \mathcal{H}$}
\smallskip

The Dirichlet-to-Neumann semigroup $\{U(t) = e^{- t \Lambda_{\gamma,\partial\Omega}}\}_{t\geq0}$
on $\partial \mathcal{H}$ is defined by self-adjoint and non-negative
Dirichlet-to-Neumann generator $\Lambda_{\gamma,\partial\Omega}$ of Proposition \ref{DtoN-H-Prop}.
\begin{Pro}\label{DtoN-H-semigr-prop-PROP}
The Dirichlet-to-Neumann semigroup $\{U(t) = e^{- t \Lambda_{\gamma,\partial\Omega}}\}_{t}$ on the
Hilbert space $\partial \mathcal{H}$ is a holomorphic quasi-sectorial contraction with values in the
trace-class $\mathfrak{C}_1 (\partial \mathcal{H})$ for  ${\rm{Re}}\, (t) > 0$.
\end{Pro}
\begin{Rem}\label{DtoN-H-semigr-prop-REMARK}
The first part of the statement follows from Proposition \ref{DtoN-H-Prop}. Since the generator
$\Lambda_{\gamma,\partial\Omega}$ is self-adjoint and non-negative, the semigroup
$\{U(t)\}_{t}$ is holomorphic and quasi-sectorial contraction for ${\rm{Re}}\, (t) > 0$, see e.g.
\cite{CaZag1}, \cite{Zag1}. Compactness of the resolvent of $\Lambda_{\gamma,\partial\Omega}$ implies the
compactness of $\{U(t)\}_{t>0}$, but to prove the last part of the statement we need a supplementary
argument about asymptotic behaviour of its eigenvalues given by the Weyl asymptotic formula, Proposition
\ref{Weyl-asympt}.
\end{Rem}
\noindent This behaviour of eigenvalues implies the second part of the Proposition
\ref{DtoN-H-semigr-prop-PROP}:
\begin{Lem}\label{DtoN-Tr-norm}
The Dirichlet-to-Neumann semigroup $U(t)$ has values in the trace-class
$\mathfrak{C}_1 (\partial \mathcal{H})$  for any $t>0$.
\end{Lem}
\noindent \begin{demo} Since the Dirichlet-to-Neumann operator $\Lambda_{\gamma,\partial\Omega}$
is self-adjoint, we have to prove that
\begin{equation}\label{semi-gr-tr-norm-bis}
\|U(t)\|_1 = \sum_{k\ge 1} \rme^{-t\lambda_k} < \infty \ ,
\end{equation}
for $t>0$. Here $\|\cdot\|_1$ denotes the norm in the
trace-class $\mathfrak{C}_1 (\partial \mathcal{H})$. Then the Weyl asymptotic formula implies that
there exists bounded $M$ and function $r(k)$ such that
\begin{align*}\sum_{k\ge 1}
\rme^{-t\lambda_k}&\le  \sum_{k\ge 1} \exp \{-t[(k/c)^{\frac
1{d-1}}+r(k)] \}\\&\le \rme^{tM}\sum_{k\ge 1} \exp \{-t(k/c)^{\frac
1{d-1}} \} .
\end{align*}
Here $c:=C(\partial \Omega, \Lambda_\gamma)$ and the last sum converges for any $t>0$, which
proves the equation (\ref{semi-gr-tr-norm-bis}).
\end{demo}

\bigskip

\noindent \textbf{2.5 Example: Lax semigroups}
\smallskip

A beautiful example of explicit representation of the Dirichlet-to-Neumann semigroup (\ref{solution})
is due to Lax
\cite{Lax}, Ch.36.

Let $\gamma(x)= I$, and $\Omega = \Omega_R $, see Section 2.2. Following \cite{Lax} we define the mapping:
\begin{equation}\label{Lax-map}
K(t): v(x)\mapsto v(e^{-t/R}\ x) \ \  {\rm{for \ any}} \ \  u \in C(\Omega_R) \ ,
\end{equation}
which is a semigroup for the parameter $t\geq 0$ in the Banach space $X=C(\Omega_R)$:
\begin{equation}\label{K-semigr}
(K(\tau)K(t)v)(x)= v(e^{-\tau/R} \ e^{-t/R}\ x) = v(e^{-(\tau+t)/R} \ x) \ , \ \tau,t \geq 0 \ ,
\ x\in \Omega_R \ .
\end{equation}
\begin{Rem}\label{Lax-K}
It is clear that if  $v(x)$ is {\rm{(}}$\gamma=I${\rm{)}}-harmonic in $C(\Omega_R)$, then the function:
$x \mapsto v(e^{-t/R}\ x)$ is also harmonic. Therefore,
\begin{equation}\label{Harm-lift}
u_f(t,x):= v_f(e^{-t/R}\ x) = (K(t)L_{\partial\Omega_R} f)(x) =
(L_{\partial\Omega_R} f_t)(x)\ , \  x \in \Omega_R \ ,
\end{equation}
is the harmonic lifting of the function $f_t(\omega):= v_f(e^{-t/R}\ \omega) \ , \ \omega \in
\partial\Omega_R$,
where $v_f$ solves the problem {\rm{\textbf{(P1)}}} for $\gamma=I$. Since in the spherical coordinates
$x = (r, \theta, \varphi)$ one has:
\begin{equation*}
\partial v_f(t,x)/{\partial t}= -\partial_r v_f(e^{-t/R}r, \theta, \varphi)e^{-t/R} \ (r/R)
\end{equation*}
and
\begin{equation*}
\partial v_f(t, R, \theta, \varphi)/{\partial \nu_I} = \partial_r v_f(e^{-t/R}r, \theta, \varphi)e^{-t/R} \ ,
\end{equation*}
we get that $\partial u_f(t,\omega)/{\partial t} + \partial u_f(t,\omega)/{\partial \nu_I}=0$, i.e.
the function (\ref{Harm-lift}) is a solution of the problem {\rm{\textbf{(P2)}}}.
\end{Rem}
\noindent Hence, according to (\ref{solution}) and (\ref{Harm-lift}) the operator family:
\begin{equation}\label{Lax-semigr}
S(t): = T_{\partial\Omega_R} K(t) L_{\partial\Omega_R} \ , \ t\geq 0 \ ,
\end{equation}
defines the Dirichlet-to-Neumann semigroup corresponding to the problem \textbf{(P2)} for $\gamma(x)= I$,
and $\Omega = \Omega_R $, which is known as the \textit{Lax semigroup}. By virtue of (\ref{Harm-lift}) and
(\ref{Lax-semigr}) the action of this semigroup is known explicitly:
\begin{equation}\label{Lax-semigr-explic}
(S(t)f)(\omega) = v_f(e^{-t/R} \omega)  \ ,  \ \omega \in \partial \Omega_R \ .
\end{equation}
Notice that the semigroup relation:
\begin{equation}\label{Lax-semigr-relation}
S(\tau)S(t)= T_{\partial\Omega_R} K(\tau) L_{\partial\Omega_R} T_{\partial\Omega_R} K(t)
L_{\partial\Omega_R}=
S(\tau + t) \ ,
\end{equation}
follows from the properties of lifting and trace operators (see Remark \ref{L-T-oper-Rem}), from identity
(\ref{K-semigr}) and definition (\ref{Lax-semigr}). One finds generator
$\Lambda_{\gamma=I,\partial\Omega_R}$ of
this semigroup from the limit:
\begin{eqnarray}\label{Lax-gener}
&& 0 = \lim_{t\rightarrow 0}\sup_{\omega\in \partial\Omega_R}| \frac{1}{t} (f - S(t)f)(\omega) -
(\Lambda_{\gamma=I,\partial\Omega_R} f)(\omega)| =\\
&&  \lim_{t\rightarrow 0}\sup_{\omega\in \partial\Omega_R}| \frac{1}{t} (v_f(R, \theta, \varphi) -
v_f(e^{-t/R}\ R, \theta, \varphi))- (\Lambda_{\gamma=I,\partial\Omega_R} f)(R, \theta, \varphi)| . \nonumber
\end{eqnarray}
Then operator
\begin{equation}\label{D-to-N-Lax}
(\Lambda_{\gamma=I,\partial\Omega_R} f)(R, \theta, \varphi) = \partial_r v_f(r = R, \theta, \varphi)
\end{equation}
for any function $f$ from domain:
\begin{equation}\label{Dom-gen-Lax}
{\rm{dom}}(\Lambda_{I,\partial\Omega_R})=\{f\in \partial C(\Omega_R):
v_f\in{\rm{Ker}}(A_{I,\partial\Omega_R}) \, {\rm{and}} \, |(\partial_r v_f)\mid_{\partial\Omega_R}|
<\infty \}
\end{equation}
is identical to (\ref{dom-Lambda}) for the case: $\gamma=I$ and $\partial\Omega = \partial\Omega_R$. Therefore,
generator (\ref{D-to-N-Lax}) of the Lax semigroup is the the Dirichlet-to-Neumann operator in this particular
case of the Banach space $\partial X = C(\partial\Omega_R)$.

Similarly we can consider the Lax semigroup (\ref{Lax-semigr}) in the Hilbert space $\partial \mathcal{H}=
L^2(\partial \Omega_R, \ d\theta \, \sin\theta \, d\varphi)$. Since generator of this semigroup is
a particular case of the Dirichlet-to-Neumann operator (\ref{DtoN=R}), by (\ref{DtoN=R-spect}) and
(\ref{v_f=R})
we again obtain the corresponding action in the explicit form:
\begin{eqnarray}\label{Lax-Hilb-explic}
&&(S(t)f)(\omega) = \\
&&(e^{- t \Lambda_{I,\partial\Omega_R}} f)(\omega)) = \sum_{l=0}^{\infty}\sum_{m=-l}^{m=l}
\sum_{s=0}^{\infty} \frac{(-t)^s}{s!}\left(\frac{l}{R}\right)^s \ f_{l,m}^{(R)}  \ Y_{l,m} (\theta, \varphi)
= \nonumber \\
&& \sum_{l=0}^{\infty}\sum_{m=-l}^{m=l} (e^{-\ t/R})^l \ f_{l,m}^{(R)}  \ Y_{l,m} (\theta, \varphi) =
v_f (e^{- t/R} \omega) \ ,  \ \omega \in \partial \Omega_R \ , \nonumber
\end{eqnarray}
which coincides with (\ref{Lax-semigr-explic}).

Notice that for $t>0$ the Lax semigroups have their values in the trace-class
$\mathfrak{C}_1 (\partial \mathcal{H})$. This explicitly follows from (\ref{DtoN=R-spect}), i.e. from the fact
that the spectrum of the semigroup generator
$\sigma(\Lambda_{I,\partial\Omega_R}) := \{\lambda_{l,m}= \ l/R\}_{l=0, |m|\leq l}^{\infty}$ is discrete and
\begin{equation}\label{Tr}
{\rm{Tr}}\, S(t) = \sum_{l=0}^{\infty}(2l + 1) \ e^{- l/R} < \infty \ .
\end{equation}
The last is proven in the whole generality in Theorem \ref{DtoN-Tr-norm}.

\section{Product approximations of Dirichlet-to-Neumann semigroups}\label{Appr}
\smallskip
\noindent \textbf{3.1 Approximating family}

\smallskip

Since in contrast to the Lax semigroup ($\gamma = I$) the action of the general Dirichlet-to-Neumann
semigroup for $\gamma \neq I$ is known only implicitly (\ref{solution}), it is useful to construct converging
approximations, which are simpler for calculations and analysis.

One of them is the Emamirad-Laadnani \textit{approximation} \cite{EmLa}, which is motivated by the
explicit action
(\ref{Lax-semigr-explic}), (\ref{Lax-Hilb-explic}) of the Lax semigroup
\begin{eqnarray}\label{Lax-semigr-T-L}
&& (S(t)f)(\omega) = (T_{\partial\Omega_R} K(t) L_{\partial\Omega_R} f) (\omega)=
v_f(e^{-t/R} \omega) , \ \omega \in \partial \Omega_R \ ,\\
&& K_R(t): v(x)\mapsto v(e^{-t/R}\ x) \ \  {\rm{for \ any}} \ \
v \in C(\Omega_R)\, ({\rm{or}}\, \mathcal{H}(\Omega_R)) \ . \nonumber
\end{eqnarray}
The suggestion of \cite{EmLa} consists in substitution of the family $\{K(t)\}_{t \geq 0}$ by the
$\gamma-$\textit{deformed} operator family:
\begin{equation}\label{K-gamma-family}
K_{\gamma,R}(t): v(x)\mapsto v(e^{-(t/R) \ \gamma(x)}\ x) \ \  {\rm{for \ any}} \ \
v \in C(\Omega_R)\, ({\rm{or}}\, \mathcal{H}(\Omega_R)) \ .
\end{equation}
\begin{Def}\label{V-appr-Def}
For the ball $\Omega_R$ the Emamirad-Laadnani approximating family
$\{V_{\gamma, R}(t):= V_{\gamma,\partial \Omega_R}(t)\}_{t \geq 0}$ is defined by
\begin{equation}\label{V-appr}
(V_{\gamma, R}(t)f)(\omega): = (T_{\partial\Omega_R} K_{\gamma}(t) L_{\partial\Omega_R} f) (\omega)=
v_f(e^{-(t/R) \ \gamma(\omega)} \  \omega) , \ \omega \in \partial \Omega_R \ .
\end{equation}
\end{Def}
\begin{Rem}\label{V-appr-REM}
{\rm{(a)}} Notice that the approximating family (\ref{V-appr}) is not a semigroup:
\begin{eqnarray}\label{not-semi}
&&(V_{\gamma, R}(t) V_{\gamma, R}(s) f)(\omega) =
(T_{\partial\Omega_R} K_{\gamma}(t) L_{\partial\Omega_R} \widetilde{f}(s)) (\omega) =\\
&& v_{\widetilde{f}(s)}(e^{-(t/R) \ \gamma(\omega)} \  \omega)\neq
v_{f}(e^{-((t+s)/R) \ \gamma(\omega)} \  \omega) = (V_{\gamma, R}(t + s))f)(\omega) \ . \nonumber
\end{eqnarray}
{\rm{(b)}} This family is strongly continuous at $t=0$:
\begin{equation}\label{V-cont}
\lim_{t\searrow 0} V_{\gamma, R}(t)f = f  \ \ {\rm{for \ any}} \ \
f \in \partial X \, ({\rm{or}}\,  \partial\mathcal{H}) \ .
\end{equation}
{\rm{(c)}} By definition (\ref{V-appr}) this family has derivative at $t= +0$:
\begin{equation}\label{V-deriv}
(\partial_{t}V_{\gamma, R}(t) f)(\omega)\mid_{t=0} =
- \nu(\omega)\cdot\gamma(\omega)(\nabla v_f )(\omega) = - (\Lambda_{\gamma,\partial\Omega_R}f)(\omega) \ ,
\end{equation}
which for any $f\in {\rm{dom}}(\Lambda_{\gamma,\partial\Omega_R})$ coincides with the (minus)
Dirichlet-to-Neumann operator (\ref{Lambda}).
\end{Rem}

\smallskip

\noindent \textbf{3.2 Strong approximation of the Dirichlet-to-Neumann semigroups}

\smallskip

By virtue of Remark \ref{V-appr-REM} the Emamirad-Laadnani approximation family verifies the conditions
of the Chernoff \textit{approximation theorem} (Theorem 1.1, \cite{Che}):
\begin{Pro}\label{Chernoff}
Let $\{\Phi(s)\}_{s \ge 0}$ be a family of linear
contractions on a Banach space $\mathfrak{B}$ and let $X_0$ be the generator of a $C_0$-contraction
semigroup. Define $X(s) := s^{-1}(I - \Phi(s))$, $s > 0$. Then for $s \rightarrow +0$ the family
$\{X(s)\}_{s > 0}$ converges strongly in the resolvent sense to the operator $X_0$  if
and only if the sequence $\{\Phi(t/n)^n\}_{n \ge 1}$, $t > 0$,  converges strongly
to $e^{-tX_0}$ as $n \rightarrow \infty$, uniformly on
any compact $t$-intervals in $\mathbb{R}^{1}_{+}$.
\end{Pro}
Notice that $\{V_{\gamma, R}(t)\}_{t \geq 0}$ in the Banach space $\partial X$ is the family of contractions
because of the \textit{maximum principle} for the $\gamma$-harmonic functions $v_f$. Since the
Dirichlet-to-Neumann operator (\ref{Lambda}) is densely defined and closed, Remark \ref{V-appr-REM} (c)
implies that the family $X(s) := s^{-1}(I - V_{\gamma, R}(s))$ converges for $s \rightarrow +0$ to
$X_0 = \Lambda_{\gamma,\partial\Omega_R}$ in the strong resolvent sense.

The similar arguments are valid for the case of the Hilbert space $\partial \mathcal{H}$. By virtue of
Remark \ref{DtoN-H-Properties} the Dirichlet-to-Neumann operator $\Lambda_{\gamma,\partial\Omega}$ is
non-negative and self-adjoint. This implies again that (\ref{V-appr}) is the family of contractions in
$\partial \mathcal{H}$ and that by Remark \ref{V-appr-REM} (c) the family $X(s) := s^{-1}(I - V_{\gamma, R}(s))$
converges for $s \rightarrow +0$ to $X_0 = \Lambda_{\gamma,\partial\Omega_R}$ in the \textit{strong} resolvent sense.

Resuming the above observations we obtain the \textit{strong} approximation of the Dirichlet-to-Neumann
semigroup $U(t)$:
\begin{Cor}\label{Srong-appr}\cite{EmLa}
\begin{equation}\label{V-strong-U}
\lim_{n\to\infty}(V_{\gamma, R}(t/n))^n f= U(t)f \, ,\quad\text{for every} \ f \in \partial X
\ \text{or}  \ \partial \mathcal{H} \ ,
\end{equation}
uniformly on any compact $t$-intervals in $(0,\infty)$.
\end{Cor}

The Emamirad-Laadnani approximation theorem (Corollary \ref{Srong-appr}) has the following important
extension to more \textit{general} geometry than ball \cite{EmLa}.
\begin{Def}\label{Non-ball}
We say that a bounded smooth domain $\Omega$ in $\mathbb R^d$ has the {\rm{property of the
interior ball}}, if for any $\omega\in \partial \Omega$ there exists a tangent to $\partial\Omega$ at
$\omega$ plane $\mathcal{T}_\omega$, and such that one can construct a ball tangent to
$\mathcal{T}_\omega$ at
$\omega$,
which is totally included in $\Omega$.
\end{Def}

If $\Omega $ has this property, then with any point $\omega\in \partial\Omega$, one can associate a
\textit{unique} point $x_\omega$, which is the center of the \textit{biggest} ball $B(x_\omega,r_\omega)$
of radius
$r_\omega$ included in $\Omega$. For any $0< r \le r_\omega$, we can construct the approximating family
$V_r(t)$ related to the ball $B(x_{r,\omega},r):=\{ x\in\Omega : |x-x_{r,\omega}|\le r\}$ of radius $r$,
which is centered on the line perpendicular to $\mathcal{T}_\omega$ at the point $\omega\in \partial\Omega$,
i.e.
$x_{r,\omega}= (r/{r_\omega})x_\omega+(1-r/{r_\omega})\omega$. Then we define
\begin{equation}\label{V-r}
(V_{\gamma, r}(t)f)(\omega) := T_{\partial \Omega} \ v_f \left( x_{r,\omega} +
\rme^{-({t}/r)\gamma(\omega)}(r \ \nu_\omega)\right) \ .
\end{equation}
Here $\nu_\omega$ is the outer-normal vector at $\omega$, the function $v_f = L_{\partial \Omega} f$ is the
$\gamma$-harmonic lifting of the boundary condition $f$ on $\partial \Omega$ , and
$T_{\partial \Omega}$ is the trace operator:
\begin{equation}\label{Tr-oper}
T_{\partial \Omega}: H^1(\Omega)\ni v  \longmapsto  v \mid_{\partial \Omega} \in
H^{1/2}(\partial\Omega).
\end{equation}
\begin{Rem}\label{Gen-Omega}
Notice that:\\
{\rm{(a)}} since $\nu_\omega = (\omega-x_{r,\omega})/r$, one gets $(V_{\gamma, r}(t=0)f)(\omega) :=
(T_{\partial \Omega} \ v_f) (\omega) = f (\omega)$ ;\\
{\rm{(b)}} by virtue of (\ref{V-r}) the strong derivative at $t=0$ has the form:
\begin{equation*}
(\partial_t V_{\gamma, r}(t=0)f)(\omega)
= - \gamma(\omega)\nu_\omega \cdot (\nabla v_f)(\omega) = - (\Lambda_{\gamma,\partial\Omega}f)(\omega) ,
\end{equation*}
see (\ref{V-deriv}).
\end{Rem}
\begin{Pro}\label{Srong-appr-Gen-Omega}\cite{EmLa}
Let $\Omega$ has the property of interior ball, and let
\begin{eqnarray*}\label{balls}
&&\inf_{\omega\in\partial\Omega} \{r>0: B(x_\omega,r_\omega)\subset \Omega \} >0 , \\
&&\sup_{\omega\in\partial\Omega} \{r>0: B(x_\omega,r_\omega)\subset \Omega \} < \infty \ . \nonumber
\end{eqnarray*}
For any $0<s\leq 1$ we define $V_{\gamma, s r_\omega}$, i.e.
\begin{equation}\label{s-family}
V_{\gamma, s r_\omega}f(\omega)=
v_f \left( x_{s,\omega} +
\rme^{-({t}/(s r_\omega))\gamma(\omega)}(s r_\omega \ \nu_\omega)\right) \ ,
\end{equation}
where $x_{s,\omega}= s x_{\omega} +(1-s)\omega$. Then for any $0<s\leq 1$
\begin{equation}\label{Lim-s}
\lim_{n\to\infty}(V_{\gamma, s r_\omega}(t/n))^n f= U(t)f \, ,\quad\text{for every} \ f \in \partial X
\ \text{or}  \ \partial \mathcal{H} \ ,
\end{equation}
uniformly on any compact $t$-intervals in $(0,\infty)$.
\end{Pro}
\begin{Rem}\label{Gen-Omega-V}
By Definition \ref{V-appr-Def} for the ball $\Omega_R$ and constant matrix-valued function $\gamma(x) = I$
one obviously have $V_{\gamma= I, R}(t) = S(t) = U(t)$. On the other hand, for a general smooth domain
$\Omega$ with geometry verifying the conditions of  Proposition \ref{Srong-appr-Gen-Omega}, one
is obliged to consider the family of approximations $V_{\gamma, s r_{\omega}}$ even for the homogeneous
case $\gamma = I$.
\end{Rem}

\section{Dirichlet-to-Neumann Gibbs semigroups}\label{Gibbs}

\smallskip

\noindent \textbf{4.1 Gibbs semigroups}

\smallskip

Since by Lemma \ref{DtoN-Tr-norm} for any  Dirichlet-to-Neumann semigroup we obtain:
$U(t>0)\in \mathfrak{C}_1 (\partial \mathcal{H})$, then one can check that it is in fact a
\textit{Gibbs} semigroup. To this end we recall main definitions and some results that we need for
the proof, see e.g. \cite{Zag2}.

Let $\mathfrak{H}$ be a separable, infinite-dimensional complex Hilbert space.
We denote by $\mathcal{L}(\mathfrak{H})$ the algebra of
all bounded operators on $\mathfrak{H}$ and by $\mathfrak{C}_\infty(\mathfrak{H})\subset\mathcal{L}(\mathfrak{H})$
the subspace of all \textit{compact} operators. The $\mathfrak{C}_\infty(\mathfrak{H})$ is a
$\ast$-\textit{ideal} in $\mathcal{L}(\mathfrak{H})$, that is: if $A \in\mathfrak{C}_\infty(\mathfrak{H})$, then
$A^* \in \mathfrak{C}_\infty(\mathfrak{H})$ and, if $A \in \mathfrak{C}_\infty(\mathfrak{H})$ and
$B \in \mathcal{L}(\mathfrak{H})$, then $A B \in \mathfrak{C}_\infty(\mathfrak{H})$ and $B A \in
\mathfrak{C}_\infty(\mathfrak{H})$. We say that a compact operator $A \in \mathfrak{C}_\infty(\mathfrak{H})$
belongs to the \textit{von Neumann-Schatten} $\ast$-ideal $\mathfrak{C}_p (\mathfrak{H})$ for a certain
$1 \le p < \infty$, if the norm
\begin{equation}\label{Jp}
\|A\|_p := \left(\sum_{n\ge 1} s_n(A)^p\right)^{1/p} < \infty ,
\end{equation}
where $s_n(A):=\sqrt{\lambda_n (A^*A)}$ are the \textit{singular} values of $A$, defined by the eigenvalues
$\{\lambda_n (\cdot)\}_{n\ge 1}$ of non-negative self-adjoint operator $A^*A$.
Since the norm $\|A\|_p$ is a non-increasing function of $p>0$, one gets:
\begin{equation}\label{_|p}
\|A\|_1 \ge \|A\|_p \ge \|A\|_q > \|A\|_{\infty} (=\|A\|) ,
\end{equation}
for $1 \leq p \leq q < \infty$. Then for the von Neumann-Schatten ideals this implies inclusions:
\begin{equation}\label{Jp_sub}
\mathfrak{C}_1 (\mathfrak{H}) \subseteq \mathfrak{C}_p (\mathfrak{H})\subseteq
\mathfrak{C}_q (\mathfrak{H}) \subset \mathfrak{C}_\infty(\mathfrak{H}) .
\end{equation}

Let $p^{-1}=q^{-1}+r^{-1}$. Then by virtue of the \textit{H\"older inequality} applied to (\ref{Jp}) one gets:
$\|AB\|_p \le \|A\|_q\|B\|_r$, if $A\in \mathfrak{C}_q (\mathfrak{H})$ and $\ B\in \mathfrak{C}_r (\mathfrak{H})$.
Consequently we obtain:
\begin{Lem}\label{ABinTr} The operator $A$ belongs to the trace-class $\mathfrak{C}_1 (\mathfrak{H})$ if and only
if there exists two {\rm{(Hilbert-Schmidt)}} operators $K_1 , \ K_2 \in  \mathfrak{C}_2 (\mathfrak{H})$,
such that $A = K_1 \ K_2$. Similarly, if $K\in  \mathfrak{C}_p (\mathfrak{H})$, then
$K^p \in  \mathfrak{C}_1 (\mathfrak{H})$.

Let $K$ be integral operator in the Hilbert space
$L^2(D, \mu)$. It is a Hilbert-Schmidt operator if and only
if its kernel $k(x,y)\in L^2(D\times D, \mu \times \mu)$ and then one gets the estimate:
$\|K\|_{2} \leq \|k\|_{L^2(D\times D, \mu \times \mu)}$.
\end{Lem}
The proof is quite straightforward and can be found in, e.g., \cite{Kat}, \cite{Sim}.
\begin{Def}\label{Gibbs-semi}\cite{Zag2}
Let $\{G(t)\}_{t\ge 0}$ be a $C_0$-semigroup on $\mathfrak{H}$ with
$\{G(t)\}_{t > 0}\subset\mathfrak{C}_{\infty}(\mathfrak{H})$. It
is called {\textit{immediate}} Gibbs semigroup, if $G(t) \in \mathfrak{C}_1 (\mathfrak{H})$
for any $t > 0$ and it is called \textit{eventually} Gibbs semigroup, if there is $t_0 > 0$, such that
$G(t) \in \mathfrak{C}_1 (\mathfrak{H})$ for any $t\ge t_0 $.
\end{Def}
\begin{Rem}\label{Gibbs-semi-REM}
{\rm{(a)}} Notice that by Lemma \ref{ABinTr} any $C_0$-semigroup such that one has
$\{G(t)\}_{t > 0}\subset\mathfrak{C}_{p}(\mathfrak{H}) $  for some $p < \infty$, is
an immediate Gibbs semigroup.\\
{\rm{(b)}} Since compact $C_0$-semigroups are norm-continuous for any $t>0$, the immediate Gibbs semigroups
are $\|\cdot\|_{1}$-norm continuous for $t>0$.
\end{Rem}
For more details of the Gibbs semigroups properties we refer to the book \cite{Zag2}.

\begin{Cor}\label{Gibbs-semi-DtoN}
By virtue of Proposition \ref{DtoN-H-semigr-prop-PROP}, Definition \ref{Gibbs-semi} and
Remark \ref{Gibbs-semi-REM} the Dirichlet-to-Neumann semigroup
$\{U(t) = e^{- t \Lambda_{\gamma,\partial\Omega}}\}_{t}$ on the Hilbert space
$\partial \mathcal{H}$ is a $\|\cdot\|_1$-holomorphic quasi-sectorial immediate Gibbs
for $\mathrm{Re}\,(t)>0$.
\end{Cor}

\smallskip

\noindent \textbf{4.2 Compact and Tr-norm approximating family}

\smallskip

\begin{Pro}\label{V-compact-Banach}\cite{EmLa}
For the ball $\Omega_R$ the Emamirad-Laadnani approximating family $\{V_{\gamma, R}(t)\}_{t \geq 0}$ consists
of compact operators on the Banach space $\partial X = C(\partial\Omega_R)$ for any $t>0$.
\end{Pro}
The proof follows from Definition \ref{V-appr-Def} by \textit{Arzela-Ascoli} criterium of compactness, since
representation (\ref{V-appr}) and conditions on $\gamma$ imply the uniform bound and equicontinuity
of the sets $\{V_{\gamma, R}(t)(\partial X )\}_{t}$ for any $t>0$.

For the case of the Hilbert space we recall the following useful condition for characterization of the
Tr-class operators \cite{Zag2}.
\begin{Pro}\label{Abstr-Tr-class}
If $A\in \mathcal{L}(\mathfrak{H})$ and $\sum_{j=1}^{\infty}\|A e_j\| < \infty$ for an orthonormal basis
$\{e_j\}_{j=1}^{\infty}$ of $\mathfrak{H}$, then $A \in \mathfrak{C}_{1}(\mathfrak{H})$.
\end{Pro}
\begin{Theo}\label{V-Tr-class}
On the Hilbert space $\partial \mathcal{H}= L^2 (\partial \Omega_R)$ the approximating family
$\{V_{\gamma, R}(t)\}_{t > 0} \subset \mathfrak{C}_1 (\partial \mathcal{H})$.
\end{Theo}
\noindent \textit{Proof}: Since the eigenfunctions $\{\phi_k\}_{k=1}^{\infty}$ of the self-adjoint
Dirichlet-to-Neumann operator $\Lambda_{\gamma,\partial\Omega_R}$ form an orthonormal basis
in $L^2 (\partial \Omega_R)$, we apply Proposition \ref{Abstr-Tr-class} for this basis.

Let $\partial \Omega_{t,\gamma,R}:=\{x_\omega := e^{-(t/R) \ \gamma(\omega)} \
\omega\}_{\omega\in\partial \Omega_R}$.
By representation (\ref{V-appr}) and  by estimate (\ref{estim}) one obtains
\begin{eqnarray}
&&\|V_{\gamma, R}(t)\phi_k \|^2= \int_{\partial \Omega_R} d\sigma(\omega)
|v_{\phi_k}(x_\omega)|^2  \nonumber \\
&&\leq  |\partial \Omega_R| \ \sup_{\omega\in\partial \Omega_R}
\psi(x_\omega,p,\partial \Omega_{t,\gamma,R})^2/k ^{2p/(d-1)} \ . \label{norm-on-basis}
\end{eqnarray}
Then by hypothesis \textbf{(H2)} on the matrix $\gamma$ one gets for the norm of the vector
$x_\omega$ in $\mathbb{R}^d$ the estimate:
\begin{equation*}
\|x_\omega\|\leq \|e^{-(t/R) \ \gamma}\| \ R \leq e^{- c_1(t/R)} \, R  \ .
\end{equation*}
Hence, for any  $t > 0$ the $\mathrm{dist}(x_\omega, \partial \Omega_R) \geq (1 - e^{- c_1(t/R)}) R > 0$,
which implies for the estimates in (\ref{estim})  and in (\ref{norm-on-basis}) that
\begin{equation*}
0<\inf_{\omega\in\partial \Omega_R}\psi(x_\omega,p,\partial \Omega_{t>0 ,\gamma, R}) \leq
\sup_{\omega\in\partial \Omega_R}
\psi(x_\omega,p,\partial \Omega_{t>0 ,\gamma, R}).
\end{equation*}
Then for $2p/(d-1) > 1$ the estimate (\ref{norm-on-basis}) ensures the convergence of the series in
the inequality:
\begin{equation*}
\|V_{\gamma, R}(t)\|_1 \leq \sum_{k=1}^{\infty}\|V_{\gamma, R}(t)\phi_k \| ,
\end{equation*}
which finishes the proof. \hfill$\Box$
\section{Concluding remarks: trace-norm approximations}\label{Concl}

The strong Emamirad-Laadnani approximation theorem (Corollary \ref{Srong-appr}) and the results of
Section 4.2 proving that Dirichlet-to-Neumann semigroup $U(t)$ and
approximants $V_{\gamma, \partial\Omega}(t/n)^n$ belong to $\mathfrak{C_1}(\partial\mathcal{H})$,
for all $n \geq 1$ and $t>0$, motivate the following conjecture:
\begin{Conj}\label{Tr-norm-conv} \cite{EmZa}
The Emamirad-Laadnani approximation theorem is valid in the $\mathrm{Tr}$-norm topology of
$\mathfrak{C_1}(\partial\mathcal{H})$.
\end{Conj}
\begin{Rem}\label{Tr-norm-conv-Rem}
Notice that the \textit{strong} approximation of the Dirichlet-to-Neumann Gibbs semigroup $U(t)$ by
the $\mathrm{Tr}$-class family $(V_{\gamma, \partial\Omega}(t/n))^n$ does not lift automatically the topology of
convergence to, e.g., operator-norm approximation \cite{Zag2}.
\end{Rem}

Therefore, to prove the Conjecture \ref{Tr-norm-conv} one needs additional arguments similar to those of
\cite{CaZag2}. To this end we put the difference in question
$\Delta_n (t):=(V_{\gamma, \partial\Omega}(t/n))^n - U(t)$ in the following form:
\begin{eqnarray}\label{Delta}
\Delta_n (t) &=& \{(V_{\gamma,\partial\Omega}(t/n))^{k_n} - (U(t/n))^{k_n}\}(V_{\gamma, R}(t/n))^{m_n} \\
 &+& (U(t/n))^{k_n} \{(V_{\gamma, \partial\Omega}(t/n))^{m_n} - (U(t/n))^{m_n}\}. \nonumber
\end{eqnarray}
Here for any $n>1$, we define two variables $k_n=[n/2]$ and
$m_n=[(n+1)/2]$, where $[x]$ denotes the \textit{integer part} of
$x\geq0$, i.e.,  $n=k_n+m_n$.
Then for the estimate of $\Delta_n (t)$ in the $\mathfrak{C_1}(\partial\mathcal{H})-$topology one gets:
\begin{eqnarray}\label{Delta-estim}
\|\Delta_n (t)\|_{1} &\leq& \|(V_{\gamma, \partial\Omega}(t/n))^{k_n} -
(U(t/n))^{k_n}\| \ \|(V_{\gamma, \partial\Omega}(t/n))^{m_n}\|_{1} \\
&+& \|(U(t/n))^{k_n}\|_{1}  \ \|(V_{\gamma, \partial\Omega}(t/n))^{m_n} - (U(t/n))^{m_n}\|. \nonumber
\end{eqnarray}

In spite of Remark \ref{Tr-norm-conv-Rem}, the \textit{explicit} representation of approximants
$\left\{(V_{\gamma, \partial\Omega}(t/n))^n\right\}_{n \geq 1}$ allows to prove the corresponding
operator-norm estimate.
\begin{Theo} \label{norm-conv} \cite{EmZa} Let  $V_{\gamma, \partial \Omega_R}(t)$ be defined by (\ref{V-appr}).
Then one gets the estimate:
\begin{equation}\label{norm-est}
\left\|(V_{\gamma, \partial \Omega_R}(t/n))^n - U(t)\right\| \leq \varepsilon(n) \ , \ \
\lim_{n\rightarrow\infty}\varepsilon(n)= 0 ,
\end{equation}
uniformly for any $t$-compact in $\mathbb{R}^{1}_{+}$.
\end{Theo}

To establish (\ref{norm-est}) we use the "telescopic" representation:
\begin{eqnarray}\label{telescop}
&&(V_{\gamma, \partial \Omega_R}(t/n))^n - U(t) = \\
&& \sum_{s=0}^{n-1} (V_{\gamma, \partial \Omega_R}(t/n))^{(n-s-1)} \{V_{\gamma, \partial \Omega_R}(t/n) - U(t/n)\}
(U(t/n))^{s} \ , \nonumber
\end{eqnarray}
and the operator-norm estimate of $\{V_{\gamma, \partial \Omega_R}(t/n) - U(t/n)\}$ for large $n$.

The next auxiliary result establishes a relation between family of operators $V_{\gamma, \partial \Omega_R}(t)$
and the Dirichlet-to-Neumann semigroup $U(t)$.
\begin{Lem}\label{V&U} \cite{EmZa} There exists a bounded operator $W_{\gamma, \partial \Omega_R}(t)$ on
$\partial \mathcal{H}$ such that
\begin{equation}\label{VWU}
V_{\gamma, \partial \Omega_R}(t)=W_{\gamma, \partial \Omega_R}(t)U(t) \ ,
\end{equation}
for any $t \geq 0$.
\end{Lem}

Now we return to the main inequality (\ref{Delta-estim}). To estimate the \textit{first} term in the
right-hand side of (\ref{Delta-estim}) we need Theorem \ref{norm-conv} and the \textit{Ginibre-Gruber}
inequality \cite{CaZag2}:
\begin{equation*}
\|(V_{\gamma, \partial\Omega}(t/n))^{m_n}\|_{1} \leq C \ U(m_n t/n) \ .
\end{equation*}
To establish the latter we use representation (\ref{VWU}) given by Lemma \ref{V&U}.

To estimate the \textit{second} term one needs only the result of Theorem \ref{norm-conv}. All together
this gives a proof of Conjecture \ref{Tr-norm-conv} at least for the ball $\Omega_R$.

\vskip 2.0cm

\noindent \textbf{Acknowledgements}

\bigskip

This paper is based on the lecture given by the author on the \textit{Lyapunov Memorial Conference},
June 24-30, 2007 (Kharkov University, Ukraine). I would like to express my gratitude to organizers and in
particular to Prof. Leonid A.Pastur for invitation and for support.

The conference talk, as well as the present account, are a part of the common project with Prof. Hassan Emamirad,
I would like to thank him for a fruitful and pleasant collaboration.

\newpage

\end{document}